\documentclass{amsart}

\newtheorem{theorem}{Theorem}[section]
\newtheorem{lemma}[theorem]{Lemma}

\theoremstyle{definition}
\newtheorem{definition}[theorem]{Definition}

\newtheorem*{note}{Notation}
\theoremstyle{remark}
\newtheorem{remark}[theorem]{Remark}

\numberwithin{equation}{section}

\newcommand{\ra}{\rightarrow}

\newcommand{\PZ}{\ensuremath{P_{Z} (\varphi)}}

\newcommand{\VPZ}{\ensuremath{P^{\ast}_{Z} (\varphi)}}

\newcommand{\PZo}{\ensuremath{P_Z ^{\mathit{classic}} (\varphi)}}

\newcommand{\ua}{\underline{a}}

\newcommand \ul {\underline}
\newcommand \bl {\begin{lemma}}
\newcommand \el {\end{lemma}}
\newcommand \bt {\begin{theorem}}
\newcommand \et {\end{theorem}} 
\newcommand \mc{\mathcal}
\newcommand \limn {\lim_{n \ra \infty}}
\newcommand \limk {\lim_{k \ra \infty}}

\newcommand \limnk {\lim_{n_k \ra \infty}}
\newcommand\IR{\mathbb R}
\newcommand \IN{\mathbb N}

\newcommand \ZZ{\mathbb Z}
\newcommand \CU{\mathcal U}

\newcommand \IG{\mathcal G}
\newcommand \BU{\mathbf U}

\newcommand \XU{\mathbf{X(U)}}
\newcommand \SU{\mathcal{S(U)}}

\newcommand \IP{P}

\newcommand \Vx{ \mathcal V (x) }
\newcommand \htop{h_{top}}
\newcommand \bp{\begin{proof}}
\newcommand \ep{\end{proof}}
\newcommand \norm {\| \varphi \|_{\infty}}

\begin{document}
\title[Topological Pressure for Non-Compact Sets]{A thermodynamic definition of topological pressure for non-compact sets}%{Irregular sets for maps with the almost specification property and for the $\beta$-transformation}
\author{Daniel Thompson}
\address{Department of Mathematics, Pennsylvania State University, University Park, State College, PA 16802, USA}
\email{thompson@math.psu.edu}
%\date{16th February 2009}
\begin{abstract}
We give a new definition of topological pressure for arbitrary (non-compact, non-invariant) Borel subsets of metric spaces. This new quantity is defined via a suitable variational principle, leading to an alternative definition of an equilibrium state. We study the properties of this new quantity and compare it with existing notions of topological pressure. We are particularly interested in the situation when the ambient metric space is assumed to be compact. We motivate the naturality of our definition by applying it to some interesting examples, including the level sets of the pointwise Lyapunov exponent for the Manneville-Pomeau family of maps.
\end{abstract}
\maketitle
\section{Introduction.}
It has long been thought desirable to generalise the standard theory of topological pressure and equilibrium states to non-compact spaces \cite{Bo2}, \cite{DJ}, \cite{GS}, \cite{HKR}, \cite{Hass}, \cite{Sa}. The purpose of this paper is to contribute a new and elementary approach to this problem. %There are currently a number of definitions of topological pressure in the non-compact setting. 

First, let us review some of the most well known and frequently used of the definitions which currently exist in the literature. Notably, Sarig developed the theory of Gurevic pressure for countable state shifts \cite{Sa}. This theory has many applications, particularly in the study of non-uniformly hyperbolic systems, where inducing schemes lead naturally to the study of countable state shifts \cite{PeSe}, \cite{BrTo}.

In another direction, Bowen defined topological entropy for non-compact subsets of a compact metric space as a characteristic of dimension type \cite{Bo2}. The study of topological entropy for the multifractal decomposition of Birkhoff averages (for example) is a well accepted goal in its own right \cite{Pe}, \cite{Bad}. Pesin and Pitskel contributed a definition of topological pressure for non-compact sets \cite{PP2}  which generalises the Bowen definition and is also suitable for the study of multifractal analysis \cite{Tho6}, \cite{Tho4}. 

Another approach is to use a definition involving the minimum cardinality of spanning sets which resembles the usual definition of topological pressure in the compact invariant setting. There are two distinct quantities which can be defined this way, called in \cite{Pe} the upper and lower capacity topological pressure (see \S \ref{1.1}).

It would be desirable for the non-compact topological pressure to satisfy a variational principle analogous to the celebrated theorem of Walters \cite[Theorem 9.10]{Wa} 
\[
\IP(\varphi) = \sup \left\{ h_{\mu} + \int \varphi d\mu \right\}.
\]
While there is a good theory of thermodynamic formalism for countable state shifts, initiated by Sarig, the techniques do not generalize to the study of arbitrary non-compact topological dynamical systems. A variational principle for the pressure of Pesin and Pitskel does exist but only applies to sets satisfying a certain condition which is very difficult to check (see \S \ref{1}). No general variational principle is known in the non-compact or non-invariant case for the upper or lower capacity topological pressure (although the relativised variational principle of Ledrappier and Walters involves the consideration of upper capacity topological pressure, see \S \ref{1.1}). In contrast, our new notion of pressure satisfies a suitable variational principle by its very definition.

The alternative notion of topological pressure for non-compact spaces which we present in this paper has an elementary definition, made via a suitable variational principle. The new pressure has the advantage that its properties are significantly easier to derive than that of the dimension-like version and we seem to pay no price in terms of desirable properties. It is of particular interest that we arrive naturally at a new definition of equilibrium state. We provide an example, described below, where this appears to be the correct notion of equilibrium state. 

While we hope that the new topological pressure will have useful applications in the future, our focus for the moment is an intrinsic study of the definition. We derive the basic properties of the new definition, we give a thorough comparison with other existing definitions, and we give some first examples which we hope will be illuminating. In particular, we give a simple example which illustrates the difference in the thermodynamic properties of the new pressure and the Pesin and Pitskel pressure. We also describe an example taken from the multifractal analysis of the Lyapunov exponent for the Manneville-Pomeau family of maps, which seems particularly well adapted to our new framework. 

We note that the new definition works best when we consider the restriction of a continuous map on a compact ambient space to a non-compact or non-invariant subset. This is the setting that we focus on for the majority of the paper.  Nevertheless,  the definition still makes sense formally when there is no ambient space or the ambient space is non-compact, and we explore this possibility in \S \ref{5}.

In \S \ref{1}, we state our definition and set up our notation. In \S \ref{2}, we study the properties of our new topological pressure when the ambient space is compact. In \S \ref{3}, we study the relationship between the different definitions. In \S \ref{4}, we  consider some interesting examples. In \S \ref{5}, we study our new topological pressure when the ambient space is non-compact. In the appendix, we prove a result which we use in \S \ref 2.
\section{Definitions and Notation.} \label{1}
Let $(X,d)$ be a compact metric space and $f:X \mapsto X$ a continuous map. Let $C(X)$ be the space of continuous real-valued functions on $X$.  Let $Z \subset X$ be an arbitrary Borel set. Let $\mathcal{M}_f (X)$ denote the space of $f$-invariant Borel probability measures on $X$ %and let $\mathcal{M}^e_f (X)$ consist of the ergodic measures. 
and $\mathcal{M}^e_f (X)$ denote those which are ergodic. %the subspace of $\mathcal{M}_f (X)$ for which the measures are ergodic. 
If $Z$ is $f$-invariant, let $\mathcal{M}_f (Z)$ denote the subset of $\mathcal{M}_f (X)$ for which the measures $\mu$ satisfy the additional property $\mu (Z) = 1$. Let $\mathcal{M}_f^e (Z) := \mathcal{M}_f (Z) \cap \mathcal{M}_f^e (X)$. For $\varphi \in C(X)$, $x \in X$ and $n \geq 1$, we write $S_n \varphi (x) := \sum_{i=0}^{n-1} \varphi (f^i x)$ and define the probability measure
\[
\delta_{x, n} = \frac{1}{n} \sum_{k=0}^{n-1} \delta_{f^k (x)}
\]
where $\delta_x$ denotes the Dirac $\delta$-measure supported on $x$. 
%Note that $\int \varphi d \delta_{x, n} = \frac{1}{n} S_n \varphi (x)$. 
We define $\mathcal{V} (x)$ to be the set of limit points for $\delta_{x, n}$, namely:
\[
\mathcal{V} (x) = \{ \mu \in \mathcal{M}_{f} (X) : \delta_{x, {n_k}} \ra \mu \mbox{ for some } {n_k} \ra \infty \}.
\]
%We sometimes use the notation $\mc V (Z) := \bigcup_{x \in Z} \Vx$. 
We state the new definition which will be the object of our study. 
\begin{definition}
Let $Z$ be an arbitrary non-empty Borel set and $\varphi \in C (X)$. Define
\[
\VPZ = \sup \left\{ h_{\mu} + \int_X \varphi d\mu \mbox{ : } \mu \in \mathcal{V} (x) \mbox{ for some } x \in Z \right\}.
\]
We set $\IP^{\ast}_\emptyset (\varphi) = \inf_{x \in X} \varphi (x)$. 
If $\varphi \equiv 0$, than we may denote $\IP^{\ast}_{Z} (0)$ by $ \htop^\ast (Z)$.
\end{definition}
\begin{note}
We denote the topological pressure of $\varphi$ on $Z$ defined as a dimension characteristic using the definition of Pesin (see \S \ref{3}) by $\PZ$ and $\htop (Z) := \IP_Z (0)$. The new topological pressure of definition 2.1 and quantities associated with it will always carry an asterisk, eg. $\VPZ$, $\htop^\ast (Z)$.
\end{note}
\begin{remark}
An alternative natural definition to make is as follows:
\[
P_{Z}^{\#} (\varphi) = \sup \left\{ h_{\mu} + \int_X \varphi d\mu \mbox{ : } \mu = \limn \delta_{x,n} \mbox{ for some } x \in Z \right\}.
\]
If no such measures exist, then we set $P^{\#}_Z (\varphi) = \inf_{x \in X} \varphi (x)$. One obvious relationship is $\VPZ \geq P_{Z}^\# (\varphi)$. We take the point of view that $\VPZ$ is the better quantity to study because it captures more information about $Z$ than $P^\#_Z (\varphi)$. Furthermore, the relationship between $\VPZ$ and $\PZ$ is better than the relationship between $P_{Z}^\# (\varphi)$ and $\PZ$ (see \S \ref{3}). Theorem \ref{irr} gives an example of a set $Z$ for which $\htop^\ast (Z) = \htop (Z) = \htop (f)$ but $P_{Z}^{\#} (0) = 0$.
\end{remark}
\begin{remark}
When the ambient space $X$ is non-compact, we can define $\htop^\ast (Z)$ as in definition 2.1, although we must insist that if $\bigcup_{x \in Z} \Vx = \emptyset$, then $\htop^\ast (Z) = 0$. The definition of $\IP_{Z}^\ast (\varphi)$ requires a small modification in the non-compact setting and we study this situation further in \S \ref{5}.
\end{remark}
We recall the variational principle for $\PZ$ proved by Pesin and Pitskel.
\begin{theorem}[Pesin and Pitskel] \label{rem}
Let $Z$ be $f$-invariant and $\mathcal{L} (Z) = \{ x \in Z : \mathcal{V} (x) \cap \mathcal{M}_f (Z) \neq \emptyset \}$. Then $\IP_{\mathcal{L} (Z)} (\varphi) = \sup \left\{ h_{\mu} + \int_Z \varphi d\mu \right \}$, where the supremum is taken over either $\mathcal{M}_f (Z)$ or $\mathcal{M}_f^e (Z)$. 
\end{theorem}
\subsection{Some classical notions.} \label{1.1}
We fix notation on some classical results which we use repeatedly. For an invariant measure $\mu$, let $G_\mu$ denote its set of generic points
\[
G_\mu = \{ x \in X : \delta_{x, n} \ra \mu \}.
\] 
If $\mu$ is ergodic, $G_\mu$ is non-empty and by Birkhoff's theorem $\mu (G_\mu) = 1$. Furthermore, if $f$ satisfies definition 2.2 (specification), $G_\mu$ is non-empty for any invariant measure. A proof using a slightly stronger specification property is available in \cite{De}, although the result holds true under definition 2.2. When $h_\mu > 0$, it is a corollary of the result $\htop (G_\mu) = h_\mu$ for any invariant measure. This was proved under weak assumptions which cover our setting in \cite{PfS}. 
\begin{definition}
A continuous map $f: X \mapsto X$ satisfies the specification property if for all $\epsilon > 0$, there exists an integer $m = m(\epsilon )$ such that for any collection $\left \{ I_j = [a_j, b_j ] \subset \IN : j = 1, \ldots, k \right \}$ of finite intervals with $dist( I_i, I_j) \geq m(\epsilon ) \mbox{ for } i \neq j$ and any $x_1, \ldots, x_k$ in $X$, there exists a point $x \in X$ such that
\begin{equation} \label{3a.02}
d(f^{p + a_j}x, f^p x_j) < \epsilon \mbox{ for all } p = 0, \ldots, b_j - a_j \mbox{ and every } j = 1, \ldots, k.
\end{equation}
\end{definition}
For a compact, invariant set $X$, we denote the classical topological pressure, defined as in \cite{Wa}, by $\IP_X^{classic} (\varphi)$. We use the notation $\htop (f) := \IP_X^{classic} (0)$. It is well known that $\IP_X^{classic} (\varphi) = \sup \left\{ h_{\mu} + \int \varphi d\mu \right\}$, where the supremum can be taken over all measures in $\mathcal{M}_f (X)$ or just the ergodic ones. We refer to this result as the classical variational principle. %The pressure defined for compact non-invariant subsets as in \cite{LW} coincides with the upper capacity topological pressure as defined in \cite{Pe}. Example 11.1 of \cite{Pe} shows that in general topological pressure and upper capacity topological pressure differ.
\begin{remark}
The usual definition of $\IP_X^{classic} (\varphi)$ in terms of spanning sets generalises to non-compact and non-invariant sets of a compact metric space. Let \[
Q_n (Z, \varphi, \epsilon) = \inf \{\sum_{x \in S} \exp S_n \varphi (x) : \mbox{ S is an $(n, \epsilon)$ spanning set for Z }\}.
\]
$\overline{CP}_Z (\varphi)$ is defined to be $\lim_{\epsilon \ra 0} \limsup_{n \ra \infty} \frac{1}{n} \log Q_n (Z, \varphi, \epsilon)$ and called in \cite{Pe} the upper capacity topological pressure. The lower capacity topological pressure $\ul{CP}_Z (\varphi)$ is given by repacing the $\limsup$ with $\liminf$. In $\S 11$ of \cite{Pe}, Pesin shows that these quantities can be formulated as characteristics of dimension type and example 11.1 of \cite{Pe} shows that they do not always coincide with $\PZ$, even for compact non-invariant sets. We note that in the context of fibred systems (i.e. $(X_1, f_1)$ and $(X_2, f_2)$ are dynamical systems and $\pi: X_1 \mapsto X_2$ continuous satisfies $\pi (X_1) = X_2$ and $\pi \circ f_1 = f_2 \circ \pi$.), the relativized variational principle of Ledrappier and Walters \cite{LW} involves the pressure of compact non-invariant sets (the fibres), and they use $\overline{CP}_Z (\varphi)$ rather than $\PZ$. We state the entropy version of the relativized variational principle: given $\nu \in \mc M_f (X_2)$,
\[
\sup_{\mu : \mu \circ \pi^{-1} = \nu} h_\mu = h_\nu + \int_{X_2} \overline{CP}_{\pi^{-1} (x)} (0) d \nu (x).
\]
%The same is true of Bowen's formula for fibred systems \cite{Bo3}, which in our notation reads $\htop (f_1) \leq \htop (f_2) + \sup_{x \in X_2} \overline{CP}_{\pi^{-1} (x)} (0)$, where 
\end{remark}
\section{Properties of $\VPZ$.} \label{2}
\bt \label{aa}
The topological pressure of definition 2.1 satisfies:

(1) $\IP^{\ast}_{Z_1} (\varphi) \leq \IP^{\ast}_{Z_2} (\varphi)$ if $Z_1 \subset Z_2 \subset X$,

%(2) $\IP^{\ast}_Z (\varphi) = \sup_{i \geq 1} \IP^{\ast}_{Z_i} (\varphi)$ where $Z = \bigcup_{i \geq 1} Z_i$,

(2) $\IP^{\ast}_Z (\varphi) = \sup \{\IP^{\ast}_{Y} (\varphi): Y \in \mc F \}$ where $Z = \bigcup_{Y \in \mc F} Y$ and $\mc F$ is a collection (countable or uncountable) of Borel subsets of $X$,

(3) $\IP^{\ast}_Z (\varphi \circ f) = \VPZ$,

(4) If $\psi$ is cohomologous to $\varphi$, then $\IP^{\ast}_Z (\varphi) = \IP^{\ast}_Z (\psi)$,

(5) $\IP^{\ast}_Z (\varphi + \psi) \leq \VPZ + \beta (\psi)$, where $\beta (\psi) = \sup_{\mu \in \mathcal{M}_f (X)} \int_X \psi d\mu$,

(6) $\IP^{\ast}_Z ((1-t) \varphi + t \psi) \leq (1-t) \VPZ + t \IP^{\ast}_Z (\psi)$.

(7) $| \IP^{\ast}_Z (\varphi) - \IP^{\ast}_Z (\psi) | \leq \| \psi - \varphi \|_\infty$,

(8) $\VPZ \geq \inf_{x \in X} \varphi(x)$,

(9) For every $k \in \ZZ$, $\IP^{\ast}_{f^k Z} (\varphi) = \VPZ$,

(10) $\IP^{\ast}_{\bigcup_{k \in \ZZ} f^k Z} (\varphi) = \IP^{\ast}_{\bigcup_{k \in \IN} f^{-k} Z} (\varphi) = \IP^{\ast}_{\bigcup_{k \in \IN} f^{k} Z} (\varphi) = \VPZ.$
\et
\bp
Since $\bigcup_{x \in Z_1} \Vx \subseteq \bigcup_{x \in Z_2} \Vx$, the first statement is immediate. The second statement is true because $\bigcup_{x \in Z} \Vx \subseteq \bigcup_{Y \in \mc F} \bigcup_{x \in Y} \Vx$. 
%The second statement follows from the fact that $\bigcup_{x \in Z} \Vx \subseteq \bigcup_i \bigcup_{x \in Z_i} \Vx$.
It is a standard result that $\Vx \subseteq \mathcal{M}_f (X)$ (see for example \cite{Wa}) and thus $\int_X \varphi d\mu = \int_X \varphi \circ f d\mu$ for $\mu \in \Vx$. The third statement follows. If $\psi$ is cohomologous to $\varphi$, then there exists a continuous function $h$ so $\psi = \varphi + h - h \circ f$ and so $\int_X \varphi d\mu = \int_X \psi d\mu$. The fourth statement follows. We leave (5) and (6) as easy exercises. %For (5), we consider the following inequality, taking supremums over $\bigcup_{x \in Z} \Vx$:
%$\sup \{ h_\mu + \int_X \varphi d\mu + \int_X \psi d\mu \} \leq \sup \{ h_\mu + \int_X \varphi d\mu \} + \sup \{ \int_X \psi d\mu \}$.
%For (6), we note that $\IP^{\ast}_Z ((1-t) \varphi + t \psi) \leq \sup \{ 2 h_\mu + (1-t) \int_X \varphi d\mu + t \int_X \psi d\mu \}$.
(7) follows from the fact that for $\mu \in\Vx$,
\[
h_\mu + \int \varphi d \mu \leq h_\mu + \int \psi d \mu + \| \psi - \varphi \|_\infty.
\] 
(8) follows from the fact that $h_\mu + \int \varphi d \mu \geq \inf_{x \in X} \varphi(x)$. (9) is true because  $\Vx = \bigcup_{ \{y: y = f^k x \} } \mc V (y)$ for all $x \in Z$ and we can apply (2). (10) follows from (9) and (2). 
\ep

%\bt \label{ad.2}
%For every $k \in \ZZ$, $\IP^{\ast}_{f^k Z} (\varphi) = \VPZ$. Thus
%\[
%\IP^{\ast}_{\bigcup_{k \in \ZZ} f^k Z} (\varphi) = \IP^{\ast}_{\bigcup_{k \in \IN} f^{-k} Z} (\varphi) = \IP^{\ast}_{\bigcup_{k \in \IN} f^{k} Z} (\varphi) = \VPZ.
%\]
%\et
$\VPZ$ is a topological invariant of dynamical systems in the following sense:
\bt \label{ab}
Let $(X_i, d_i)$ be compact metric spaces and $f_i : X_i \mapsto X_i$ be continuous maps for $ i = 1,2$. Let $\pi : X_1 \mapsto X_2$ be a homeomorphism satisfying $\pi \circ f_1 = f_2 \circ \pi$. Then for any continuous $\varphi : X_2 \mapsto \IR$ and Borel $Z \subset X_2$, we have $\VPZ = \IP^{\ast}_{\pi^{-1}(Z)} (\varphi \circ \pi)$. 
\et
\bp
For $\psi \in C (X_2)$ and $\mu \in \mathcal{M}_{f_2} (X_2)$, let $\tilde \psi := \psi \circ \pi$ and $\tilde \mu := \mu \circ \pi$. 
%Note that $\tilde \psi \in C (X_1)$ and $\int \tilde \psi d \tilde \mu = \int \psi d \mu$. 
Let $\mu \in \bigcup_{x \in Z} \Vx$. Then $\mu = \lim_{n_k \ra \infty} \delta_{x, n_k}$ for some $x \in Z$, $n_k \ra \infty$. Let $y \in X_1$ satisfy $\pi (y) =x$. For an arbitrary function $\psi \in C (X_1)$, 
\begin{eqnarray*}
\int \psi d \tilde \mu & = & \int \psi \circ \pi^{-1} d \mu \\ &=& \limnk \frac{1}{n_k} S_{n_k} \psi \circ \pi^{-1} (x) \\&=& \limnk \frac{1}{n_k} S_{n_k} \psi (y)\\ &=& \limnk \int \psi d \delta_{n_k, y}. 
\end{eqnarray*}
Since this is true for all $\psi \in C(X_1)$, we have $\tilde \mu \in \mc V (y)$. Thus $\mu \in \bigcup_{x \in Z} \Vx \Rightarrow \tilde \mu \in \bigcup_{y \in \pi ^{-1} (Z)} \mc V (y)$. Since $h_{\tilde \mu} + \int \tilde \varphi d \tilde \mu =  h_{\mu} + \int \varphi d\mu$, then $\IP^{\ast}_{\pi^{-1}(Z)} (\tilde \varphi) \geq \VPZ$. Reversing the previous argument gives the desired equality.
\ep
The proof shows that if $\pi$ were only assumed to be a continuous surjective map, we would obtain the inequality $\VPZ \leq \IP^{\ast}_{\pi^{-1}(Z)} (\varphi \circ \pi)$. We now verify that in the compact, invariant case $\VPZ$ agrees with the classical topological pressure. 
\bt \label{ac}
If $Z$ is compact and $f$-invariant, then $\VPZ = \PZo$.
\et
\bp
By compactness of $Z$, $\mathcal{M}_f (Z)$ is compact and thus $\bigcup_{x \in Z} \Vx \subseteq \mathcal{M}_f (Z)$. The inequality $\VPZ \leq \PZo$ follows immediately. For the opposite inequality, let $\mu \in \mathcal{M}_f (Z)$ be ergodic. %By Birkhoff's theorem, the set $G_\mu = \{ x \in Z : \delta_{x, n} \ra \mu \}$ satisfies $\mu (G_\mu) = 1$. 
Taking any point $x$ in  $G_\mu$, we have $\Vx = \mu$. We conclude that $\mathcal{M}_f^e (Z) \subseteq \bigcup_{x \in Z} \Vx$ and the desired inequality follows from the classical variational principle.
\ep
The following result is clear from the definition.
\bt \label{ac.1}
Suppose $Z$ contains a periodic point $x$ with period $n$. Then we have $\IP^{\ast}_{\{x\}} (\varphi) = \frac{1}{n} \sum_{i=0}^{n-1} \varphi (f^i x)$ and $\VPZ \geq \frac{1}{n} \sum_{i=o}^{n-1} \varphi (f^i x)$.
\et
%We say a point is generic for an invariant measure $\mu$ if $\delta_{x, n} \ra \mu$, and denote by $G_{\mu}$ the set of all such points. One can easily verify that if $\mu$ is ergodic, then $\mu (G_{\mu}) = 1$. 
We now consider the set of generic points $G_\mu$. Bowen (for entropy \cite{Bo2}) and Pesin (for pressure \cite{PP2}) showed that $\IP_{G_\mu} (\varphi) = h_{\mu} + \int \varphi d\mu$. In fact, it was this property that motivated Bowen's original dimensional definition of topological entropy. We see that similar properties holds for the new topological pressure.
\bt \label{ac.2}
For any invariant measure, $\IP^{\ast}_{G_\mu} (\varphi) = h_{\mu} + \int \varphi d\mu$. Let $Z$ be a Borel set with $Z \cap G_\mu \neq \emptyset$, then $\VPZ \geq h_{\mu} + \int \varphi d\mu$. Now assume that $\mu$ is an equilibrium measure for $\varphi$, then $\IP^{\ast}_{G_\mu} (\varphi) = \IP_X^{classic} (\varphi)$. In particular, let $m$ be a measure of maximal entropy and $Z \cap G_m \neq \emptyset$. Then $h_{top}^{\ast} (Z) = h_{top} (f)$.
\et
The proof follows immediately from the definitions. Let us remark that if a measure of maximal entropy is fully supported then $h_{top}^{\ast} (U) = h_{top} (f)$ for every open set $U$. 

It is informative to consider the pressure of a single point.
\bt \label{ac.3}
Let $x \in G_\mu$. Then $\IP^{\ast}_{\{x\}} (\varphi) = h_{\mu} + \int \varphi d\mu$ and $\IP_{\{x\}} (\varphi) = \int \varphi d\mu$. Thus $\IP^{\ast}_{\{x\}} (\varphi) = \IP_{\{x\}} (\varphi)$ iff $h_\mu =0$. 
\et
\bp
The first statement is clear. The second follows from the formula for pressure at a point $\IP_{ \{x\} } (\varphi ) = \liminf_{n \ra \infty} \frac{1}{n} S_n \varphi (x)$ (see the appendix). Since $x \in G_\mu$, $\frac{1}{n} \sum_{i = 0}^{n-1} \varphi (f^i (x)) \ra \ \int \varphi d\mu$ for every continuous $\varphi$.
\ep
\bt \label{ac.4}
Let $x \in X$. If $h_\mu > 0$ for some $\mu \in \Vx$, then $\IP^{\ast}_{\{x\}} (\varphi) > \IP_{\{x\}} (\varphi)$. %If $h_\mu = 0$ for all $\mu \in \Vx$, then $\IP^{\ast}_{\{x\}} (\varphi) = \IP_{\{x\}} (\varphi)$. 
\et
\bp
Suppose $\mu \in \mc V (x)$. Then for some $m_k \ra \infty$, we have
\[
\int \varphi d \mu =  \limk \frac{1}{m_k} S_{m_k} \varphi (x) \geq \liminf_{n \ra \infty} \frac{1}{n} S_n \varphi (x) = \IP_{\{x\}} (\varphi).
\]
Therefore, if $h_\mu > 0$, then $\IP^{\ast}_{\{x\}} (\varphi) \geq h_{\mu} + \int \varphi d\mu > \IP_{\{x\}} (\varphi)$.% and the first statement is proved.% Now let $\nu \in \Vx$ and $h_\nu = 0$. For some $l_k \ra \infty$, $\int \varphi d\nu = \limk \frac{1}{l_k} S_{l_k} \varphi (x) \leq \limsup_{n \ra \infty} \frac{1}{n} S_n \varphi (x)$ and the second statement is proved.
\ep
\begin{remark}
Theorem \ref{ac.3} provides us with a simple example which shows that $\PZ$ and $\VPZ$ are not equal. In theorem \ref{app}, we verify that for $x \in G_\mu$, $\ul {CP}_{\{x\}} (\varphi) = \overline {CP}_{\{x\}} (\varphi) = \int \varphi d\mu$. Hence, theorem \ref{ac.3} shows that $\VPZ$ cannot be equal to these quantities either.
\end{remark}
\begin{remark}
We note that $\VPZ$ is sensitive to the addition of a single point to the set $Z$. When $\varphi \neq 0$, the same is true of $\PZ$. %(We note that $\PZ \geq \sup \{ \int \varphi d \mu : Z \cap G_\mu \neq \emptyset$\}.) 
However, in the case of entropy, we have a contrast between $\htop(Z)$, which remains the same under the addition of a countable set, and $\htop^\ast(Z)$, where a single point can carry full entropy. %To avoid this, we could alter the definition to
%\[
%\htop^{\ast \ast} (Z) = \sup \{h_\mu : \mu \in \Vx \mbox{ for uncountably many } x \in Z  \}.
%\]
\end{remark}
%The following theorem gives us an interpretation of what our definition is telling us, the proof being immediate from theorem \ref{aa}(2).
%\bt
%Let $\widehat X = \{x \in X : \delta_{x,n} \mbox{ does not converge } \}$ and $\widehat X(\varphi, t) = \{x \in \widehat X : h_\mu + \int \varphi d \mu \leq t \mbox{ for all } \mu \in \Vx\}$. Then
%\[
%\VPZ = \max \left\{ \sup\{ h_\mu + \int \varphi d \mu: Z \cap G_\mu \neq \emptyset \}, \sup\{t : Z \cap \widehat X(\varphi, t) \neq \emptyset\} \right \}.
%\]
%\et
%\begin{remark}
%Let us assume that $\VPZ \cap \widehat X = \emptyset$. Then $\VPZ$ is completely determined by the sets $G_\mu$ which have non-empty intersection with $Z$ and for which $h_\mu + \int \varphi d \mu$ is largest.
%\medbreak
%The following theorem gives us an interpretation of what our definition is telling us, the proof being immediate from 
%\bt
%Let $\widehat X = \{x \in X : \delta_{x,n} \mbox{ does not converge } \}$ and $\widehat X(\varphi, t) = \{x \in \widehat X : h_\mu + \int \varphi d \mu \leq t \mbox{ for all } \mu \in \Vx\}$. Then
%\[
%\VPZ = \max \left\{ \sup\{ h_\mu + \int \varphi d \mu: Z \cap G_\mu \neq \emptyset \}, \sup\{t : Z \cap \widehat X(\varphi, t) \neq \emptyset\} \right \}.
%\]
%\et
%\begin{remark}
%Suppose $Z \subseteq \bigcup_{\mu \in \mc M_f (X)} G_\mu$. It follows from theorem \ref{aa}(2) that $\VPZ = \sup\{ h_\mu + \int \varphi d \mu: Z \cap G_\mu \neq \emptyset \}$. Thus, $\VPZ$ is completely determined by which sets $G_\mu$ have non-empty intersection with $Z$.
%\end{remark}
For ergodic measures, an inverse variational principal holds.
\bt \label{ad}
Suppose $\mu$ is ergodic. Then

(1) $h_{\mu} = \inf \{ h_{top}^{\ast} Z : \mu (Z) = 1 \}$,

(2) $h_{\mu} + \int \varphi d\mu = \inf \{ \VPZ : \mu (Z) = 1 \}$.
\et
\bp
We prove (2), then (1) follows as a special case. Suppose $Z$ is a Borel set with $\mu (Z) = 1$. Since $\mu$ is assumed to be ergodic, $\mu (G_{\mu}) = 1$ and thus $Z \cap G_\mu \neq \emptyset$. It follows that $\VPZ \geq h_{\mu} + \int \varphi d\mu$  and thus $\inf \{ \VPZ : \mu (Z) = 1 \} \geq h_{\mu} + \int \varphi d\mu$. Since $\IP^{\ast}_{G_\mu} (\varphi) = h_{\mu} + \int \varphi d\mu$, we have an equality.
\ep
The assumption that $\mu$ is ergodic is essential. For example, let $\mu = p \mu_1 + (1-p) \mu_2$ where $\mu_1, \mu_2$ are ergodic with $h_{\mu_1} \neq h_{\mu_2}$ and $p \in (0, 1)$. If $\mu(Z) =1$, then $\mu_1 (Z) =1$ and thus $Z$ contains generic points for $\mu_1$. Therefore, $h_{top}^\ast (Z) \geq h_{\mu_1}$. Repeating the argument for $\mu_2$, we obtain $\inf \{ h_{top}^{\ast} Z : \mu (Z) = 1 \} \geq \max \{ h_{\mu_1}, h_{\mu_2} \} > h_\mu =  p h_{\mu_1} + (1-p) h_{\mu_2}$. In fact, since $\mu(G_{\mu_1} \cup G_{\mu_2}) =1$ and $\htop^\ast (G_{\mu_1} \cup G_{\mu_2}) = \max \{ h_{\mu_1}, h_{\mu_2} \}$, we have $\inf \{ h_{top}^{\ast} Z : \mu (Z) = 1 \} = \max \{ h_{\mu_1}, h_{\mu_2} \}$.

We have a version of Bowen's equation. 
\bt \label{ae}
Let $\varphi$ be a strictly negative continuous function. Let $\psi : \IR \mapsto \IR$ be given by $\psi (t) := \IP^{\ast}_Z (t \varphi)$. Then the equation $\psi (t) = 0$ has a unique solution. The solution lies in $[0, \infty)$.
\et
\bp
Let $s >t$. Let $\mu \in \bigcup_{x \in Z} \mc V (x)$ and $C = \inf -\varphi(x) > 0$.  We have
\[
h_\mu + \int s \varphi d\mu = h_\mu + \int t \varphi d\mu - (s-t) \int - \varphi d \mu
\]
and, since $\int - \varphi d \mu \in [ C, \norm]$,
\[
h_\mu + \int s \varphi d\mu \leq h_\mu + \int t \varphi d\mu - (s-t) C.
\]
Therefore, $\psi (s) - \psi(t) \leq - (s-t)C$ and so $\psi$ is strictly decreasing. (Similarly, $\psi (s) - \psi(t) \geq - (s-t) \norm$, so $\psi$ is bi-Lipschitz.) Since $\psi (0) \geq 0$, $\IP^{\ast}_Z (t \varphi) = 0$ has a unique root. 
\ep
\begin{remark}
We compare the properties derived here with those satisfied by $\PZ$. In theorem \ref{aa}, properties (1), (3), (4), (6) and (7) hold for $\PZ$. Property (2) holds for $\PZ$ only when the union is at most countable. Properties (9) and (10) are known to hold for $\PZ$ when $f$ is a homeomorphism. Theorems \ref{ab}, \ref{ac}, \ref{ad} and \ref{ae} hold for $\PZ$.
\end{remark}

\subsection{Equilibrium states for $\VPZ$.}
Suppose a measure $\mu^{\ast}$ satisfies $\VPZ = h_{\mu^{\ast}} + \int_X \varphi d\mu^{\ast}$ and $\mu^{\ast} \in \bigcup_{x \in Z} \Vx$ for a (not necessarily invariant) Borel set $Z$. Then we call  $\mu^{\ast}$ a $\ast$-equilibrium state for $\varphi$ on $Z$. If $\mu^{\ast}$ satisfies $\htop^\ast (Z) = h_{\mu^{\ast}}$, we call $\mu^{\ast}$ a measure of maximal $\ast$-entropy. If $Z$ is invariant, we call a measure $\mu$ that satisfies both $\PZ = h_{\mu} + \int_X \varphi d\mu$ and $\mu (Z) = 1$ simply an equilibrium state for $\varphi$ on $Z$. The latter definition coincides with that of Pesin \cite{Pe}. It is clear from the definition that if $\mu^{\ast}$ is a $\ast$-equilibrium state and $\mu$ is an equilibrium state for $\varphi$ on $Z$, then
\[
h_{\mu^{\ast}} + \int_X \varphi d\mu^{\ast} \geq h_{\mu} + \int_X \varphi d\mu.
\]
Note that it is possible that $\mu^{\ast} (Z) = 0$. There are situations where the new definition seems more appropriate than the old. We describe a non-trivial example in \ref{Kalpha} but first let us a consider a periodic point ${x}$ of period $n > 1$. Then, for any function, $\delta_{x,n}$ is a $\ast$-equilibrium state on $\{x \}$. However, as $\{ x \}$ is not invariant, the notion of equilibrium state is not defined.

\section{The relationship between $\PZ$ and $\VPZ$.} \label{3}
In theorem \ref{3main}, we show that the inequality $\PZ \leq \VPZ$ holds. Theorem \ref{ac.3} provides examples where $\PZ < \VPZ$ and non-trivial examples can be constructed. $\S \ref{4}$ contains concrete examples where $\PZ = \VPZ$ and we have the following: 
\bt
For an $f$-invariant Borel set $Z$, let $\mc G (Z) = \bigcup_{\mu \in \mc M_f (Z)} G_\mu \cap Z$. Then $\IP_{\mc G (Z)} (\varphi) = \IP_{\mc G (Z)}^{\ast} (\varphi)$.
\et
\bp
Note that $\mc L (G (Z)) = G (Z)$. Applying theorem \ref{rem}, we have $\IP_{\mc G (Z)} (\varphi) = \sup \{ h_\mu + \int \varphi d \mu : \mu \in \mc M_f (\mc G (Z)) \} = \IP_{\mc G (Z)}^{\ast} (\varphi)$.
\ep
Before embarking on a sketch proof that $\PZ \leq \VPZ$, we give a less sharp result, whose proof is straight forward given theorem \ref{rem}.
\bt
If $Z$ is an $f$-invariant Borel set, we have
\[
\IP_{\mathcal{L} (Z)} (\varphi) \leq \VPZ \leq  \IP_{\overline Z}^{classic} (\varphi) \mbox{ and } \IP_{\mathcal{L} (Z)} (\varphi) \leq \IP^{\ast}_{\mathcal{L} (Z)} (\varphi).
\]
\et
\bp
We note that if $\mu \in \mathcal{M}_f^e (Z)$, then $\mu (Z \cap G_\mu) = 1$. Taking $x \in Z \cap G_\mu$, we have $\Vx = \{ \mu \}$ and thus $\mathcal{M}_f^e (Z) \subseteq \bigcup_{x \in Z} \Vx $. Note that $x \in \mc L (Z)$ and so $\mathcal{M}_f^e (Z) \subseteq \bigcup_{x \in \mc L (Z)} \Vx$. By theorem \ref{rem}, the first and third inequalities follows. For the second inequality, we have $\VPZ \leq \IP_{\overline Z}^\ast (\varphi) = \IP_{\overline Z}^{classic} (\varphi)$. 
%For the third inequality, we note that $\mathcal{M}_f (Z) \subseteq \bigcup_{x \in \mc L (Z)} \Vx$ but, a priori, the inclusion and thus the inequality may be strict.
\ep
Example \ref{Kalpha} shows that the second inequality may be strict (the sets $K_\alpha$ are dense but do not carry full entropy), and the remark after lemma \ref{spec} shows that the third inequality may be strict.  The first inequality of the following theorem is the main result of this section. We do not assume that $Z$ is invariant. 
\bt \label{3main}
Let $Z$ be an arbitrary Borel set and $Y = \overline {\bigcup_{k \in \IN} f^{-k} Z}$, then
\[
\PZ \leq \VPZ \leq \IP^{classic}_{Y} (\varphi).
\]
\et
\S \ref{3.2} constitutes a sketch proof of the first inequality. This result, although never stated before, follows from part of Pesin and Pitskel's proof of theorem \ref{rem}, with only minor changes required. For a complete proof, we refer the reader to \cite{PP2} or \cite{Pe}. Here, we attempt to convey the key technical ingredients. The second inequality is trivial as $Y$ is a closed invariant set containing $Z$.
\subsection{Definition of Pesin and Pitskel's topological pressure.}
Let $(X,d)$ be a compact metric space, $f:X \mapsto X$ be a continuous map and $\varphi \in C(X)$.  Let $Z \subset X$ be a Borel subset.  We take a finite open cover $\mathcal{U}$ of X and denote by $\mathcal{S}_m (\mathcal{U})$ the set of all strings $\mathbf{U} = \{(U_{i_0}, \ldots, U_{i_{m-1}}) : U_{i_j} \in \mathcal{U} \}$ of length $m = m(\mathbf{U})$.  We define $\mathcal{S(U)} = \bigcup_{m \geq 0} \mathcal{S}_m (\mathcal{U})$, where $S_0 (\CU)$ consists of $\emptyset$. To a given string $\mathbf{U} = (U_{i_0}, \ldots, U_{i_{m-1}}) \in \mathcal{S} (\mathcal{U})$, we associate the set 
$\mathbf{X(U)} = \{ x \in X : f^j (x) \in U_{i_j} \mbox{ for all } j = 0, \ldots , m(\mathbf{U}) -1\} =  \bigcap_{j=0}^{m(\mathbf{U}) -1} f^{-j} U_{i_j}$.
We say that a collection of strings $\mathcal{G} \subset \mathcal{S(U)}$ covers $Z$ if $Z \subset \bigcup_{\mathbf{U} \in \mathcal{G}} \mathbf{X(U)}$. Let $\alpha \in \IR$.  We make the following definitions:
\[
Q(Z,\alpha, \mathcal{U}, \mathcal{G}, \varphi) = \sum_{\mathbf{U} \in \mathcal{G}} \exp \left(-\alpha m(\mathbf{U})+\sup_{x\in \mathbf{X(U)}} \sum_{k=0}^{m(\mathbf{U})-1} \varphi(f^{k}(x))\right),
\]

\[
M(Z, \alpha, \mathcal{U}, N, \varphi) = \inf_{\mathcal{G}} Q(Z,\alpha, \mathcal{U}, \mathcal{G}, \varphi),
\]
where the infimum is taken over all finite or countable subcollections of strings $\mathcal{G} \subset \mathcal{S(U)}$ such that $m(\mathbf{U}) \geq N$ for all $\mathbf{U} \in \mathcal{G}$ and $\mathcal{G}$ covers $Z$. When $\XU = \emptyset$, we set $\sup_{x\in \mathbf{X(U)}} \sum_{k=0}^{m(\mathbf{U})-1} \varphi(f^{k}(x)) = - \infty$ . Define
\[
m(Z, \alpha, \mathcal{U}, \varphi) := \lim_{N \rightarrow \infty} M(Z, \alpha, \mathcal{U}, N, \varphi).
\]
There exists a critical value $\alpha_{c}$ with $ - \infty \leq \alpha_{c} \leq + \infty$ such that $m(Z, \alpha, \mathcal{U}, \varphi) = \infty$ for $\alpha < \alpha_{c}$ and $m(Z, \alpha, \mathcal{U}, \varphi) = 0$ for $\alpha > \alpha_{c}$.
Let $| \CU | = \max \{ \mbox{Diam} U_i : U_i \in \CU \}$.
\begin{definition} 
We define the following quantities:

$(1) \IP_Z (\varphi, \CU) = \inf \{ \alpha : m(Z, \alpha, \mathcal{U}, \varphi) = 0\} = \sup \{ \alpha :m(Z, \alpha, \mathcal{U}, \varphi) = \infty \} = \alpha_c$,

$(2)\IP_Z (\varphi) : = \lim_{ | \CU |  \rightarrow 0} \IP_Z (\varphi, \CU)$.
\end{definition}

For well definedness of $\PZ$, we refer the reader to \cite{PP2} or \cite{Pe}.

\subsection{Sketch proof of $\PZ \leq \VPZ$} \label{3.2}
Let $\CU = \{U_1, \ldots, U_r\}$ be an open cover of $X$ and $\epsilon > 0$. Let 
\[
\mbox{Var} (\varphi, \CU) = \sup \{ | \varphi(x) - \varphi(y) | : x, y \in U \mbox{ for some } U \in \CU\}. 
\]
Let $E$ be a finite set of cardinality $n$, and $\underline{a} = (a_0, \ldots, a_{k-1}) \in E^k$. Define the probability vector $\mu_{\ua}= (\mu_{\ua}(e_1), \ldots, \mu_{\ua}(e_{n}))$ on $E$ by 
\[
\mu_{\ua}(e_i) = \frac{1}{k}(\mbox{the number of those $j$ for which } a_j = e_i).
\]
Define
\[
H(\ua) = - \sum_{i = 1}^{n} \mu_{\ua}(e_i) \log \mu_{\ua}(e_i).
\]
In \cite{Pe}, the contents of the following lemma are proved under the assumption that $\mu \in \mathcal{V} (x) \cap \mc M_f (Z)$. However, the property $\mu (Z) =1$ is not required. We omit the proof.

\begin{lemma} \label{lemma 3}
Given $x \in Z$ and $\mu \in \mathcal{V} (x)$, there exists a number $m > 0$ such that for any $n > 0$ one can find $N > n$ and a string $\BU \in \SU$ of length $N$ satisfying:

(1) $x \in \XU$,

(2) $\sup_{x\in \mathbf{X(U)}} \sum_{k=0}^{N-1} \varphi(f^{k}(x)) \leq N \left(\int \varphi d\mu + \mbox{Var} (\varphi, \CU) + \epsilon\right)$,

(3) $\BU = (U_0, \ldots, U_{N-1})$ contains a substring $\BU^\prime$ with the following properties: There exists $k \in \IN$ with $N-m \leq km \leq N$ and $0 \leq i_0 \leq \ldots \leq i_{k-1}$ so $a_0 = (U_{i_0}, \ldots, U_{i_0 + m})$, $\ldots, a_{k-1} = (U_{i_{k-1}}, \ldots, U_{i_{k-1} + m})$ and $\BU^\prime = (a_0, \ldots, a_{k-1})$. Note that the length of $\BU^\prime$ is $km$. Writing $E = \{a_0, \ldots, a_{k-1}\}$ and $\ua = (a_0, \ldots, a_{k-1})$, then
\[
\frac{1}{m} H(\ua) \leq h_\mu + \epsilon.
\]
\end{lemma}
Given a number $m > 0$, denote by $Z_m$ the set of points $x \in Z$ for which there exists a measure $\mu \in \mathcal{V} (x)$ so lemma \ref{lemma 3} holds for this $m$. We have that $Z = \bigcup_{m>0} Z_m$. Denote by $Z_{m, u}$ the set of points $x \in Z_m$ for which there exists $\mu \in \mathcal{V} (x)$ so lemma \ref{lemma 3} holds for this $m$ and $\int \varphi d\mu \in [u-\epsilon, u+\epsilon]$. Set $c = \sup \{ h_{\mu} + \int \varphi d\mu : \mu \in \bigcup_{x \in Z} \Vx \}$. Note that if $x \in Z_{m, u}$, then the corresponding measure $\mu$ satisfies 
\begin{equation}\label{2o}
h_{\mu} \leq c - \int \varphi d\mu \leq c- u + \epsilon.
\end{equation}
Suppose a finite set $\{u_1, \ldots, u_s\}$ forms an $\epsilon$-net of the interval $[-\| \varphi \|, \| \varphi \| ]$. Then 
\[
Z = \bigcup_{m = 1}^{\infty} \bigcup_{i = 1}^{s} Z_{m, u_i}
\]
and hence $\PZ \geq \sup_{m, i} \IP_{Z_{m, u_i}} (\varphi)$. It will suffice to prove that for arbitrary $m \in \IN$ and $u \in \IR$ that $\IP_{Z_{m, u}} (\varphi) \leq c$.

For each $x \in Z_{m, u}$, we construct a string $\BU_x$ and substring $\BU^{\prime}_x$ satisfying the conditions of lemma \ref{lemma 3}. Let $\IG_{m, u}$ denote the collection of all such strings $\BU_x$ and $\IG_{m, u}^{\ast}$ denote the collection of all such substrings $\BU^{\prime}_x$. Choose $N_0$ so $m(\BU_x) \geq N_0$ for all $\BU_x \in \IG_{m, u}$. Let $\IG_{m, u, N}$ denote the subcollection of strings $\BU_x \in \IG_{m, u}$ with $m(\BU) = N$ and $\IG_{m, u, N}^{\ast}$ denote the correponding subcollection of substrings. Note that
\[
\IG_{m, u} = \bigcup_{N=N_0}^{\infty} \IG_{m, u, N} \mbox{ and } \# \IG_{m, u, N} \leq \# \CU^m \# \IG_{m, u, N}^{\ast}.
\]
We use the following lemma of Bowen \cite{Bo}. %(There is actually a small mistake in Bowen's proof. The mistake is the omission of an error term which disappears after taking a limit. It is no big deal.)
\begin{lemma} \label{2k.1}
Fix $h > 0$. Let $R(k, h, E) = \{ \ua \in E^k : H(\ua) \leq h\}$. Then
\[
\limsup_{k \ra \infty} \frac{1}{k} \log \#(R(k, h, E)) \leq h.
\]
\end{lemma}
Set $h = c-u + \epsilon$. It follows from (\ref{2o}) and the third statement of lemma \ref{lemma 3} that if $x \in Z_{m, u}$ has an associated string $\BU_x$ of length $N$, then its substring $\BU^{\prime}_x$ is contained in $R(k, m(h + \epsilon), \CU^m)$ where $k$ satisfies $N > km \geq N-m$. Therefore, $\# \IG_{m, u, N}^{\ast}$ does not exceed $\# R(k, m(h + \epsilon), \CU^m)$, and thus $\# \IG_{m, u, N} \leq \# \CU^m \#(R(k, m(h + \epsilon), \CU^m))$. Applying lemma \ref{2k.1}, we obtain
\begin{eqnarray*}
\limsup_{N \ra \infty} \frac{1}{N} \log \# \IG_{m, u, N} & \leq & \limsup_{k \ra \infty} \frac{1}{mk} \log \#\CU^m \#(R(k, m(h + \epsilon), \CU^m))\\& \leq &h + \epsilon.
\end{eqnarray*}
Since the collection of strings $\IG_{m, u}$ covers the set $Z_{m, u}$, we use property (2) of lemma \ref{lemma 3} to get
\begin{eqnarray*}
Q(Z_{m,u}, \lambda, \mc U, \IG_{m,u}, \varphi) &=& \sum_{N=N_0}^{\infty} \sum_{\BU \in \IG_{m,u, N}} \exp \left\{-\lambda N + \sup_{x\in \mathbf{X(U)}} \sum_{k=0}^{N-1} \varphi(f^{k}(x))\right\}\\
&\leq& \sum_{N=N_0}^{\infty} \# \IG_{m, u, N} \exp\left\{N\left(-\lambda + \mbox{Var} (\varphi, \CU) + \int \varphi d\mu + \epsilon\right)\right\}.
\end{eqnarray*}
%\begin{eqnarray*}
%M(Z_{m,u}, \lambda, \varphi, \CU, N_0) &\leq& Q(Z_{m,u}, \lambda, \IG_{m,u})\\
%&=& \sum_{N=N_0}^{\infty} \sum_{\BU \in \IG_{m,u, N}} \exp \left\{-\lambda N + \sup_{x\in \mathbf{X(U)}} \sum_{k=0}^{N-1} \varphi(f^{k}(x))\right\}\\
%&\leq& \sum_{N=N_0}^{\infty} |\IG_{m, u, N}| \exp\left\{N\left(-\lambda + \gamma(\CU) + \int_Z \varphi d\mu + \epsilon\right)\right\}.
%\end{eqnarray*}
Choose $N_0$ sufficiently large so for $N \geq N_0$, we have $\# \IG_{m, u, N} \leq \exp (N(h +2\epsilon ))$ and thus
\[
M(Z_{m,u}, \lambda, \CU, N_0, \varphi) \leq \sum_{N=N_0}^{\infty} \exp\left\{N\left(h -\lambda + 
\mbox{Var} (\varphi, \CU) + \int \varphi d\mu + 3\epsilon\right)\right\}.
\]
Let $\beta = \exp\left(h -\lambda + \mbox{Var} (\varphi, \CU)  + \int \varphi d\mu + 3\epsilon\right)$. If $\lambda > c + \mbox{Var} (\varphi, \CU) + 5\epsilon$, then $0 < \beta < 1$. Thus,
\[ 
M(Z_{m,u}, \lambda, \CU, N_0, \varphi) \leq \frac{\beta^{N_0}}{1 - \beta},
\]
\[ 
m(Z_{m,u}, \lambda, \CU, \varphi) \leq \lim_{N_0 \ra \infty} \frac{\beta^{N_0}}{1 - \beta} = 0.
\]
It follows that $\lambda \geq \IP_{Z_{m, u}} (\varphi, \CU)$. Since we can choose $\lambda$ arbitrarily close to $c + \mbox{Var} (\varphi, \CU)  + 5\epsilon$, it follows that
\[
\IP_{Z_{m, u}} (\varphi, \CU) \leq c + \mbox{Var} (\varphi, \CU)  + 5\epsilon.
\]
We are free to choose $\epsilon$ arbitrarily small, so on taking the limit $| \CU | \ra 0$, we have $\IP_{Z_{m, u}} (\varphi) \leq c$, as required. It follows that $\IP_{Z} (\varphi) \leq c$.

\begin{remark}
In \cite{PP2}, it is shown that if $\mu \in \mc M_f (X)$ and $\mu (Z) = 1$ then $\PZ \geq h_\mu + \int \varphi d \mu$. Thus, if $Z$ is a set satisfying $\mu (Z) = 1$ for all $\mu \in \bigcup_{x \in Z} \Vx$, then $\PZ =\VPZ$.
\end{remark}
\begin{remark}
If $\PZ < \VPZ$, then we see a phenomenon similar to example \ref{Kalpha}, where probability measures $\mu$ with $\mu (Z) < 1$ or even $\mu (Z) = 0$ capture information about the set $Z$. This may seem unusual but example \ref{Kalpha} motivates the utility of this point of view.
\end{remark}
\begin{remark} \label{Phash}
We can adapt the proof to obtain the inequality $P_{\mc G (Z)} \leq \IP_{Z}^\# (\varphi)$. The argument would differ in the paragraph above lemma \ref{2k.1}. We would construct strings $\BU_x$ and $\BU_x^\prime$ only for those $x \in \mc G(Z)$ rather than every $x \in Z$.
\end{remark}
\begin{remark}
We can view the result of this section as an inequality for $\PZ$. We state this explicitly without reference to definition 2.1. Let $Z$ be a Borel subset (not necessarily invariant) of a compact metric space $(X, d)$. Then
\[
\PZ \leq \sup \{ h_\mu + \int \varphi d \mu : \mu = \limnk \delta_{x, n_k} \mbox { for some } x \in Z, n_k \ra \infty \}.
\]
\end{remark}

\section{Examples.} \label{4}
Here are some interesting examples for which $\PZ$ and $\VPZ$ coincide.
\subsection{North-South map.} \label{NSmap}
The following example was suggested by Pesin. Let $X = S^1$, $f$ be the North-South map and $Z = S^1 \setminus \{S\}$. (By the North-South map, we mean the map $f = g^{-1} \circ h \circ g$ where $g$ is the stereographic projection from a point $N$ onto the tangent line at $S$, where $S$ is the antipodal point of $N$, and $h : \IR \mapsto \IR$ is $h(x) = x/2$.) One can verify that if $x \in S^1 \setminus \{ N, S \}$, then $\Vx = \delta_S$ and it is clear that $\mathcal V (\{ N \} ) = \delta_N$. Using this and the fact that $h_{\delta_S} = h_{\delta_N} = 0$, we have
\[
\VPZ = \max \{ \int \varphi d \delta_S, \int \varphi d \delta_N \} = \max \{ \varphi(N), \varphi(S) \}
\]
To calculate $\PZ$, one can use $\IP_{Z} (\varphi) = \max \{\IP_{ \{ N \} } (\varphi), \IP_{Z \setminus \{N \}} (\varphi) \}$. Using the formula for pressure at a point or Pesin's variational principle, $\IP_{ \{ N \} } (\varphi) = \varphi (N)$. One can verify that $\IP_{Z \setminus {N}} (\varphi) = \varphi (S)$. Thus, $\PZ$ and $\VPZ$ coincide for all continuous $\varphi$.
\begin{remark}
Note that $\mc L (Z) = \{ N \}$. If we choose $\varphi$ so that $\varphi(S) > \varphi (N)$, we are furnished with an example where $\PZ > \IP_{\mathcal{L} (Z)} (\varphi)$, showing that we could not replace $\IP_{\mathcal{L} (Z)} (\varphi)$ by $\PZ$ in Pesin's variational principle (see theorem \ref{rem}). 
\end{remark}
\begin{remark}
Our example shows that, in contrast to the compact case, the wandering set can contribute to the pressure (whether we consider $\VPZ$ or $\PZ$). Let $\mc {NW} (X)$ be  the non-wandering set of $(X, f)$. % and $\mc  {W} (X):= X \setminus \mc {NW} (X)$. 
(Recall that $x \in \mc {NW} (X)$ if for any open set $U$ containing $x$ there exists $N$ so $f^N(U) \cap U \neq \emptyset$.) For an arbitrary set $Y \subset X $, let $\mc {NW} (Y) = Y \cap \mc {NW} (X)$. %and $\mc {W} (Y) = Y \cap \mc {W} (X)$. 
For the set $Z$ of our example, $\mc{NW} (Z) = N$ (see $\S 5.3$ of \cite{Wa}). Assuming that $\varphi(S) > \varphi (N)$, we have
\[
P^{\ast}_{\mc{NW}(Z)} (\varphi) = \varphi (N) < \varphi (S) = \VPZ.
\]
This contrasts with the compact case, where $P^{classic}_{\mc{NW}(X)} (\varphi) = P^{classic}_{X} (\varphi)$. 
\end{remark}
\subsection{Irregular Sets.}
\bt \label{irr}
Let $(\Sigma, \sigma)$ be a topologically mixing subshift of finite type and $\hat{\Sigma} $ be the set of non-typical points, namely:
\[
\hat{\Sigma} := \Sigma \setminus \bigcup_{\mu \in \mathcal{M}_f (\Sigma)} G_\mu.
\]
Then $h_{top}^{\ast} ( \hat{\Sigma} )= h_{top} (\sigma)$ and $ \IP_{\hat{\Sigma}}^\ast (\psi) = \IP_{\Sigma}^{classic} (\psi)$ for all $\psi \in C(X)$.
\et
We remark that Barreira and Schmeling showed in \cite{BS} that $h_{top} (\hat{\Sigma}) = h_{top}( \sigma)$. It follows that $h_{top} ( \hat{\Sigma} ) = h_{top}^{\ast} ( \hat{\Sigma} )$. After an application of the classical variational principle, the proof of theorem \ref {irr} follows immediately from the next lemma in which, for simplicity, we assume $\Sigma$ is a full shift.
\bl \label{shift}
$\mathcal{M}_f ^e(\Sigma) \subseteq \bigcup_{x \in \hat{\Sigma}} \Vx.$
\el
\bp
Let $\mu_1$ be some ergodic measure. Let $\mu_2$ be some other ergodic measure. Let $x \in G_{\mu_1}$, $y \in G_{\mu_2}$ and $N_k \ra \infty$ sufficiently rapidly that $N_{k+1} > 2^{N_k}$. We can use the specification property of the shift to construct a point $p$ so $\delta_{p, N_{2k}} \ra \mu_1$ and $\delta_{p, N_{2k+1}} \ra \mu_2$. Namely, let $w_{2i -1} = (x_1, \ldots, x_{N_{2i-1}})$ and $w_{2i} = (y_1, \ldots, y_{N_{2i}})$ for all $i \geq 1$. Let $p = w_1 w_2 w_3 \ldots \in \Sigma$. Then $p \in \hat{\Sigma}$ and $\mu_1 \in \mc V (p)$.
\ep
An analogous result holds in a more general setting.
\bt \label{irr2}
Let $(X,d)$ be a compact metric space and $f:X \mapsto X$ be a continuous map with the specification property. Let $\varphi: X \mapsto \IR$ be a continuous function satisfying $\inf_{\mu \in \mathcal{M}_{f} (X)} \int \varphi d \mu < \sup_{\mu \in \mathcal{M}_{f} (X)} \int \varphi d \mu$.
%Assume $\varphi \in C(X)$ satisfies $\inf \int \varphi d \mu < \sup \int \varphi d \mu$, where the $\inf$ and $\sup$ are over $\mathcal{M}_{f} (X)$.
%\[
%\inf_{\mu \in \mathcal{M}_{f} (X)} \int \varphi d \mu < \sup_{\mu \in \mathcal{M}_{f} (X)} \int \varphi d \mu.
%\]
Let
\[
\widehat X_{\varphi, f} := \left \{ x \in X : \lim_{n \ra \infty} \frac{1}{n} \sum_{i = 0}^{n-1} \varphi (f^i (x)) \mbox{ does not exist } \right \}.
\]
Then $h_{top}^{\ast} ( \widehat X_{\varphi, f} )= h_{top}(f)$ and $P_{\widehat X_{\varphi, f}}^\ast (\psi) = P_{X}^{classic} (\psi)$ for all $\psi \in C(X)$.
\et
%It is shown in \cite{EKL}, and independently in \cite{Tho3}, that $h_{top} ( \widehat X_{\varphi} )= h_{top} (f)$ and so we have $h_{top}^{\ast} ( \widehat X_{\varphi} ) = h_{top} ( \widehat X_{\varphi} )$. 
Under the same assumptions, it is shown in \cite{EKL} that $h_{top} ( \widehat X_{\varphi, f} )= h_{top} (f)$ and in \cite{Tho4} that $P_{\widehat X_{\varphi, f}}(\psi) = P_{X}^{classic} (\psi)$. Thus, we have $P_{\widehat X_{\varphi, f} }(\psi) = P_{\widehat X_{\varphi, f} }^\ast(\psi)$. The proof of theorem \ref {irr2} follows immediately from the next lemma by the classical variational principle.
\bl \label{spec}
$\mathcal{M}_f ^e(X) \subseteq \bigcup_{x \in \widehat X_{\varphi, f}} \Vx.$
\el
\bp[Sketch Proof] 
Let $\mu_1, \mu_2$ be ergodic measures with $\int \varphi d \mu_1 < \int \varphi d \mu_2$. Let $x_i$ satisfy $\frac{1}{n} S_n \varphi (x_i) \ra \int \varphi d \mu_i$ for $i=1,2$. Let $m_k := m(\epsilon / 2^k)$ be as in the definition of specification and $N_k$ be a sequence of integers chosen to grow to $\infty$ sufficiently rapidly that $N_{k+1} > \exp \{\sum_{i=1}^k (N_i + m_i)\}$. We define a sequence of points $z_i \in X$ inductively using the specification property. Let $d_n(x,y) = \max \{d (f^i x, f^i y) : i = 0, \ldots, n-1 \}$. Let $t_1 = N_1$, $t_{k} = t_{k-1} + m_k + N_{k}$ for $k \geq 2$ and $s(k) := (k+1) (\mbox{mod} 2) +1$. Let $z_1 = x_1$. Let $z_2$ satisfy $d_{N_1}(z_2, z_1) < \epsilon /4$ and $d_{N_2} (f^{N_1 +m_2}z_2 , x_2) < \epsilon / 4$. Let $z_k$ satisfy
$d_{t_{k-1}}(z_{k-1}, z_k) < \epsilon / 2^k$ and $d_{N_k} (f^{t_{k-1} +m_k}z_k, x_{s(k)}) < \epsilon / 2^k$.
%\[
%d_{t_{k-1}}(z_{k-1}, z_k) < \epsilon / 2^k \mbox{ and } d_{N_k} (f^{t_{k-1} +m_k}, x_{s(k)}) < \epsilon / 2^k.
%\]
Let $B_n (x, \epsilon) = \{ y \in X : d_n (x, y) < \epsilon \}$. We can verify that $\overline B_{t_{k+1}} (z_{k+1}, \epsilon /2^k) \subset \overline B_{t_{k}} (z_k, \epsilon /2^{k-1})$. Define $p := \bigcap \overline B_{t_{k}} (z_k, \epsilon /2^{k-1})$.
For any $\psi \in C(X)$, we can show $\frac{1}{t_k} S_{t_k} \psi (p) \ra \int \psi d \mu_{s(k)}$. Thus $\delta_{p, t_{2k-1}} \ra \mu_1$, $\delta_{p, t_{2k}} \ra \mu_2$ and so $\mu_1, \mu_2 \in \mc V (p)$. In particular, $p \in \widehat X_{\varphi, f}$.
\ep
\begin{remark}
Using a similar construction to the proof of lemma \ref{shift}, we can show that the inequality $\IP_{\mathcal{L} (Z)} (\varphi) \leq \IP^{\ast}_{\mathcal{L} (Z)} (\varphi)$ may be strict. Let $(\Sigma, \sigma)$ be a Bernoulli shift. Let $\mu_1, \mu_2$ be ergodic measures with $h_{\mu_1} > h_{\mu_2}$. We can construct a point $z$ so the sequence of measures $\delta_{z, n}$ does not converge and $\mc V (z) = \{ \mu_1, \mu_2 \}$. Let $Z = G_{\mu_2} \cup \{z\}$. We see that $\mc L (Z) = Z$ and, by theorem \ref{rem}, $\htop (Z) = h_{\mu_2}$. However, $\htop^{\ast} (Z) = h_{\mu_1}$. 
\end{remark}
In \cite{PfS}, Pfister and Sullivan consider a weak specification property which they call the $g$-almost product property. The construction of lemma \ref{spec} generalises unproblematically to this setting and thus the statement of theorem \ref{irr2} holds for continuous maps with the $g$-almost product property. From \cite{PfS}, we know that the $\beta$-shift satisfies the $g$-almost product property for every $\beta > 1$  (while the specification property only holds for a set of $\beta$ of zero Lebesgue measure \cite{Bu}). As a consequence of this discussion, we obtain
\bt
For $\beta > 1$, let $(X, \sigma_\beta)$ be the $\beta$-shift. Assume $\varphi \in C(X)$ satisfies $\inf_{\mu \in \mathcal{M}_{\sigma_\beta} (X)} \int \varphi d \mu < \sup_{\mu \in \mathcal{M}_{\sigma_\beta} (X)} \int \varphi d \mu$. Then $h_{top}^{\ast} ( \widehat X_{\varphi, \sigma_\beta} )= \log \beta$.
\et
\subsection{Levels sets of the Birkhoff Average.} 
\bt \label{ba}
Let $(X,d)$ be a compact metric space, $f:X \mapsto X$ be a continuous map with the specification property and $\varphi, \psi \in C (X)$. For $\alpha \in \IR$, let \[ K_\alpha = \left \{ x \in X : \lim_{n \ra \infty} \frac{1}{n} \sum_{i = 0}^{n-1} \varphi (f^i (x)) = \alpha \right \}.\]
Suppose $K_\alpha \neq \emptyset$, then 

(1) $\htop^{\ast} (K_\alpha) = \sup \left \{ h_\mu : \mu \in \mathcal M_f (X) \mbox{ and } \int \varphi d\mu = \alpha \right \}$,

(2) $ \IP_{K_\alpha}^\ast (\psi) = \sup \left \{ h_\mu + \int \psi d \mu: \mu \in \mathcal M_f (X) \mbox{ and } \int \varphi d\mu = \alpha \right \}$.
\et
In \cite{Tho6}, we studied $P_{K_\alpha} (\psi)$ under the same assumptions as theorem \ref{ba}, proving an analogous result. This shows that $P_{K_\alpha} (\psi) = P^\ast_{K_\alpha} (\psi)$ when $f$ has specification. % $h_{top} ( K_\alpha ) = h_{top}^{\ast} ( K_\alpha)$. We also note that, in a symbolic dynamics setting, Luzia \cite{Lu} proved a result analogous to (2). 
The proof of theorem \ref{ba} follows from the next lemma.
\bl \label{113}
$\{ \mu \in \mathcal M_f (X) : \int \varphi d\mu = \alpha \} = \{ \mu \in \Vx : x \in K_\alpha \}$.
\el
\bp
Let $\mu \in \mathcal M_f (X)$ and $\int \varphi d\mu = \alpha$. Recall that $G_\mu \neq \emptyset$ and let $x \in G_\mu$. Then $\Vx = \mu$, and so $\{ \mu \in \mathcal M_f (X) : \int \varphi d\mu = \alpha \} \subseteq \{ \mu \in \Vx : x \in K_\alpha \}$. Conversely, if $\mu \in \Vx$ for $x \in K_\alpha$ then there exists $n_k \ra \infty$ so $\int \varphi d\mu = \lim_{n_k \ra \infty} \int \varphi d \delta_{x, n_k} = \lim_{n_k \ra \infty} \frac{1}{n_k} S_{n_k} \varphi (x) = \limn \frac{1}{n} S_{n} \varphi (x) = \alpha$.
\ep
\subsection{Manneville-Pomeau Maps.} \label{Kalpha}
Manneville-Pomeau maps are the family of maps on $[0, 1]$ given by
\[ 
f_s (x) = x + x^{1+s} \mod 1
\]
where $s \in (0,1)$ is a fixed parameter value. Each of these maps is a topological factor of a full one-sided shift on 2 symbols and so satisfies the specification property. Takens and Verbitskiy have performed a multifractal analysis for the function $\varphi (x) = \log f_s^{\prime} (x)$ (i.e. the multifractal analysis of pointwise Lyapunov exponents). We recall some results which can be found in \cite{TV}. Let $K_\alpha$ be as before. One of the key results used for their multifractal analysis, restated in our new language, is
\bt \label{inf}
$f:X \mapsto X$ be a continuous map with the specification property, and $\varphi: X \mapsto \IR$ a continuous function. Then
\[
(1)\mbox{ }  h_{top}^{\ast} ( K_\alpha) \leq \inf_{q \in \IR} \{ \IP^{classic}_X (q \varphi) - q \alpha \}.
\]
Furthermore, if $f$ has upper semi-continuous entropy map then
\[
(2)\mbox{ }  h_{top}^{\ast} ( K_\alpha) = \inf_{q \in \IR} \{\IP^{classic}_X (q \varphi) - q \alpha\}.
\]
\et
Since $f_s$ is positively expansive, it has upper semi-continuous entropy map. There is an interval of values $\mc I$ (which turns out to be $(0, h_{\mu})$ where $\mu$ is the absolutely continuous invariant measure for $f_s$) which has the following property. For $\alpha \in \mc I$, the infimum of theorem \ref{inf} (2) is attained uniquely at $q = -1$ and $\IP^{classic}_X (- \varphi) = 0$ (using results from \cite{U} and \cite{PS}). Thus, $h_{top}^{\ast} ( K_\alpha) = \alpha$ and if $\nu$ is an equilibrium measure for $-\varphi$ with $\int \varphi d \nu = \alpha$, then $h_{top}^{\ast} ( K_\alpha) = h_\nu$. The set $\mc A = \{p\delta_0 + (1-p) \mu: p \in [0,1] \}$ consists of equilibrium measures for $-\varphi$ and Takens and Verbitsky show there is a unique measure satisfying $\mu_\alpha \in \mc A$ and $\int \varphi d \mu_\alpha = \alpha$. By lemma \ref{113}, $\mu_\alpha \in \mc V (x) $ for some $x \in K_\alpha$, and so $\mu_\alpha$ is a $\ast$-equilibrium measure (for $0$ on $K_\alpha$). However, even though $h_{top} ( K_\alpha) = h_{\mu_\alpha}$, they show $\mu_\alpha (K_\alpha) = 0$, so  $\mu_\alpha$ is not an equilibrium measure (for $0$ on $K_\alpha$) under the definition of Pesin.

In fact, $\mu_\alpha$ is the unique $\ast$-equilibrium measure. 
In Proposition 1 of \cite{PSY}, Pollicott, Sharp and Yuri show that $\nu$ is an equilibrium state for $- \varphi$ iff $\nu \in \mc A$ (they also give a nice proof that $\IP^{classic}_X (- \varphi) = 0$). It follows that if $\mu \notin \mc A$ and $\int \varphi d \mu = \alpha$, then $h_\mu < \alpha$. Combining this with the above discussion shows that $\mu_\alpha$ is unique.

\section{Topological pressure in a non-compact ambient space.} \label{5}
We define $\VPZ$ for an arbitrary set $Z \subset X $ and $\varphi \in C(X)$ when the ambient space $X$ is non-compact. For the definition to make sense, we must exclude the consideration of measures $\mu$ such that both $h_\mu = \infty$ and $\int \varphi d \mu = - \infty$. 
\begin{definition}
Let $Z$ be an arbitrary Borel set and $\varphi \in C (X)$. Define
\[
\VPZ = \sup \left\{ h_{\mu} + \int_X \varphi d\mu \mbox{ : } \mu \in \bigcup_{x \in Z} \Vx \mbox{ and } \int_X \varphi d\mu > - \infty \right\}.
\]
If $\bigcup_{x \in Z} \Vx = \emptyset$, let $\IP_{Z}^\ast (\varphi) = \inf_{x \in X} \varphi (x)$. We set $\IP^{\ast}_Z (\varphi) = -\infty$ if $\bigcup_{x \in Z} \Vx \neq \emptyset$ and $\{ \mu \in \bigcup_{x \in Z} \Vx : \int_X \varphi d\mu > - \infty \} = \emptyset$.
\end{definition}
The reason we set $\IP_{Z}^\ast (\varphi) = \inf_{x \in X} \varphi (x)$ when $\bigcup_{x \in Z} \Vx = \emptyset$ is to ensure that the inequality
$P_{Z_1}^{\ast} (\varphi) \leq P_{Z_2}^{\ast} (\varphi)$ holds for all $Z_1 \subseteq Z_2$.
We remark that if $\varphi$ is bounded below, then we have $\int_X \varphi d\mu > - \infty$ for all $\mu \in \mc M_f (X)$. Hence, if $X$ is compact, definitions 6.1 and 2.1 agree.
\begin{remark}
Assume $\htop^\ast (Z) < \infty.$  Then we do not have to restrict ourselves to measures with $\int_X \varphi d\mu > -\infty$ in the definition of $\VPZ$. Either $\VPZ = - \infty$ or the extra measures considered do not contribute to the supremum. %Let $Z$ be an arbitrary Borel set and $\varphi \in C (X)$. If $\int_X \varphi d\mu = - \infty$ for all $\mu \in \bigcup_{x \in Z} \Vx$, then $\VPZ = - \infty$. Otherwise, 
\end{remark}
%\bt
%Suppose $X$ is $\sigma$- finite. Then $\VPZ = \sup \{ \IP^{\ast}_{Z \cap K, K} (\varphi) : K \subset X \mbox{ compact, invariant } \}$.
%\et
\begin{remark}
In the non-compact setting, dimensional definitions of pressure have the disadvantage that there are examples of metrizable spaces $X$ (eg. countable state shifts) and metrics $d_1, d_2$ on $X$ where $\IP_{Z, X_1} (\varphi) \neq \IP_{Z, X_2} (\varphi)$ (where $X_1 = (X, d_1)$ and $X_2 = (X, d_2)$) but $d_1$ and $d_2$ give rise to the same Borel structure on $X$ and thus no variational principle can hold. However, $\VPZ$ depends only on the Borel structure of $X$ and is thus invariant under a change of topologically equivalent metric.
\end{remark}
\begin{remark}
In \cite{DJ}, Dai and Jiang study a definition of topological entropy for non-compact spaces adapted to the problem of estimating the Hausdorff dimension of the space. Their definition is not a topological invariant, so is not equivalent to ours. They give an interesting discussion of the issues one faces when considering entropy as a measure of chaotic behaviour in the non-compact setting.
%They remark that for the map $x \mapsto x +1$ on $\IR$, $\htop (\IR) = \infty$. They argue that this is a disadvantage 
\end{remark}
We now study some properties of $\VPZ$ in the non-compact setting. 
\bt
Let $\IP^{\ast}_{Z, Y} (\varphi)$ denote the pressure of $\varphi$ on $Z$ when $Z \subset Y$ and $Y$ is considered as the ambient space in the definition. Let $K \subset X$ be compact and invariant and $Z \subset K$. Then $\IP^{\ast}_{Z, X} (\varphi) = \IP^{\ast}_{Z, K} (\varphi)$.
\et
\bp
It suffices to notice that if $\mu \in \bigcup_{x \in Z} \Vx$, then $\mu \in \mc M_f (K)$ and $h_\mu (f |_K) = h_\mu (f)$.
\ep
\bt \label{nc}
Let $X$ be a separable metric space and $\varphi \in C(X)$. Then 

(1) $\IP_X^\ast (\varphi) = \sup \{ \IP_{K, X}^\ast (\varphi) : K \subset X \mbox{ is compact}\}$,

(2) $\IP_X^\ast (\varphi) = \sup \{ h_\mu + \int \varphi d \mu : \mu \in \mc M_f (X), \int \varphi d \mu > - \infty \}$.
\et
\bp
For (1), we note that if $K_n$ is a countable collection of compact sets that cover $X$, then $\IP_X^\ast (\varphi) = \sup \{ \IP_{K_n, X}^\ast (\varphi) \}$ by basic properties of $\IP_X^\ast (\varphi)$. For (2), let $c$ denote the value taken by the supremum. That $\IP_X^\ast (\varphi) \leq c$ is immediate. It suffices to consider only ergodic measures in the supremum. We note that since $X$ is a separable space, if $\mu$ is ergodic then $\mu (G_\mu) = 1$. Thus, there exists $x$ satisfying $\Vx = \mu$, which shows that $\IP_X^\ast (\varphi) \geq c$.
\ep

In \cite{GS}, Gurevich and Savchenko study two definitions of topological pressure adapted to non-compact spaces. We compare these with $\VPZ$.
\begin{definition}
Set $\IP^{int}(X, \varphi) = \sup \{ \IP^{classic}_K (\varphi) \}$, where the supremum is over all subsets $K \subset X$ which are compact and invariant.  Suppose $X$ can be continuously embedded in a compact metric space $\hat X$ and that $f$ and $\varphi$ can be extended continuously to $\hat X$. We set $\IP^{ext} (X, \varphi) = \inf \{ \IP_{X, \hat X} (\varphi) \}$, where the infimum is over all such embeddings.
\end{definition}
\begin{remark}
We can make a definition analogous to $P^{ext} (X, \varphi)$ but using $P^\ast_{X, \hat X} (\varphi)$ in place of $P_{X, \hat X} (\varphi)$. We denote this quantity by $P^{ext \ast} (X, \varphi)$.
\end{remark}
\bt \label{ie}
For any $X$ separable, $f: X \mapsto X$ and $\varphi \in C(X)$, we have $\IP^{int}(X, \varphi) \leq \IP_X^\ast (\varphi)$. When  $\IP^{ext} (X, \varphi)$ and $P^{ext \ast} (X, \varphi)$ are well defined, $\IP_X^\ast (\varphi) \leq \IP^{ext}(X, \varphi) \leq P^{ext \ast} (X, \varphi)$.
\et
\bp
The first inequality follows from the classical variational principle and (2) of theorem \ref{nc}. %For $\IP^{ext} (X, \varphi)$ to be well defined, we require that $\varphi$ is bounded. 
Let $\hat X$ be a compact metric space satisfying the requirements of the definition amd $\hat \varphi$ be the extension by continuity of $\varphi$ to $\hat X$. By theorem \ref{rem} and (2) of theorem \ref{nc},
\[
\IP_X^\ast (\varphi) = \IP_{\mc L (X), \hat X} (\hat \varphi) \leq \IP_{X, \hat X} (\hat \varphi).
\]
Since $\hat X$ was arbitrary, we obtain the second inequality. The third inequality is a consequence of theorem \ref{3main}.
\ep
\begin{remark}
Both inequalities of theorem \ref{ie} may be strict. As noted in \cite{GS} and \cite{HKR}, let $Y$ be a compact metric space and $f:Y \mapsto Y$ be a minimal homeomorphism with $\htop (f) > 0$. Let $\varphi = 0$. Let $X = Y \setminus \mc O(x)$, where $\mc O(x)$ is the orbit of an arbitrary $x \in Y$. There are no compact, invariant, non-empty subsets of $X$,  so $\IP^{int}(X, 0) = 0$. However, $\htop^\ast (X) = \sup \{ h_\mu : \mu \in \mc M_f (Y) \} = \htop (f)$. For the second inequality, we use an example similar to \ref{NSmap}. Let $X = S^1 \setminus\{S\}$ with induced metric $d$ from $S^1$, $f$ is the North-South map and $\varphi (x) = d (x, \{N\})$. We have $\IP_X^\ast (\varphi) = \varphi (N) = 0$. We can verify that given any continuous embedding into $\hat X$ and any $y \in \hat X \setminus X$, $\IP_{X, \hat X} (\varphi) \geq \varphi (y) > 0$. 
\end{remark}
\begin{remark}
In \cite{Hass}, the authors compare various definitions of topological entropy for a non-compact space $X$ and a continuous map $f:X \mapsto X$. One of these definitions is a natural generalisation of Adler-Konheim-McAndrew's original definition of entropy \cite{AKM}, which we denote by $\htop^{AKM} (f)$. Proposition 5.1 of \cite{Hass} provides an example of a homeomorphism $f$ of the open unit interval (equipped with a non-standard metric) for which $\htop^{AKM} (f) = \infty$ but $\htop^\ast (f) = 0$.
\end{remark}

In \cite{HKR}, Handel, Kitchens and Rudolph give another definition of entropy for a non-compact metric space $(X, d)$ and a homeomorphism $f: X \mapsto X$, which is invariant under a change of topologically equivalent metric and is a generalisation of $\overline {CP}_Z (0)$. Let $S(K, n, \epsilon, d)$ denote the smallest cardinality of an $(n, \epsilon)$ spanning set for a compact set $K \subset X$ in the metric $d$. Let
\[
\htop^{d} (X) := \sup \{\lim_{\epsilon \ra 0} \limsup_{n \ra \infty} \frac{1}{n} \log S(K, n, \epsilon, d) : K \subset X \mbox{ is compact } \}.
\]
In fact, this definition first appeared in \cite{Bo3}. The innovation of \cite{HKR} is to define
\[
\htop^{HKR} (X) := \inf \{ \htop^{d^\prime} (X) : d^\prime \mbox{ is a metric topologically equivalent to d} \}.
\]
%and in our notation reads
%\[
%\htop^{hkr} (X) = \sup\{ \overline {CP}_K (0) : K \subset X \mbox{ is compact }\}.
%\]
They show that $\htop^{HKR} (X) \geq \sup \{h_\mu : \mu \in \mc M_f (X) \}$ and construct an example where the inequality is strict. Thus $\htop^{HKR} (X) \geq \htop^\ast (X)$ and it is possible that the two quantities may not coincide. 

%However, if $X$ is locally compact, $f: X \mapsto X$ is uniformly continuous, $Y$ is the one point compactification of $X$ and $g: Y \mapsto Y$ is the extension by continuity of $f$, they show that $\htop^{HKR} (X) = \htop^\ast (X) = \htop (g)$.

\subsection{Countable state shifts of finite type.}
We conclude by considering a topologically mixing countable state shift of finite type $(\Sigma, \sigma)$. Following Sarig \cite{Sa}, we equip $\Sigma$ with the metric $d(x, y) = r^{t(x,y)}$ where $t(x,y) = \inf(\{k : x_k \neq y_k\}\cup \infty)$ and $r \in (0,1)$. Let $\IP^G (\varphi)$ denote the Gurevic pressure as defined by Sarig \cite{Sa} where $\varphi$ is a locally H\"older function and $h^G (\sigma) := \IP^G (0)$. In \cite{GS}, the authors allow $\Sigma$ to be equipped with more general metrics and study $\IP^{int}(\Sigma, \varphi)$ and $\IP^{ext}(\Sigma, \varphi)$ for $\varphi \in C(\Sigma)$. %One would like to say '$\IP_G (\varphi) = \IP_\Sigma^\ast (\varphi)''$ in full generality but the situation requires some care and it is crucial that we use the metric $d$. 
To rephrase corollary 1 of \cite{Sa}, Sarig showed that in his setting $\IP^G (\varphi) = \IP^{int}(\Sigma, \varphi)$.

\bt
$\htop^{\ast} (\Sigma) = h^G (\sigma)$.
\et
\bp
By corollary 1.7 of \cite{GS}, $\IP^{int}(\Sigma, 0) = \IP^{ext}(\Sigma, 0)$ in the metric $d$. The result follows from theorem \ref{ie}. 
\ep
\bt
We have $\IP_\Sigma^\ast (\varphi) \geq \IP^G (\varphi)$ and thus if $\IP^G (\varphi) = \infty$, then $\IP_\Sigma^\ast (\varphi) = \infty$. Under the extra assumption $\sup_{x \in \Sigma} | \sum_{\sigma y = x}  e^{\varphi (y) }| < \infty$, we have $\IP_\Sigma^\ast (\varphi) = \IP^G (\varphi) < \infty$.
\et
\bp
The first inequality is a rephrasing of theorem \ref{ie}. Under the extra assumption, Sarig showed $\IP^G (\varphi) = \sup \{ h_\mu + \int \varphi d \mu : \mu \in \mc M_\sigma (\Sigma), \int \varphi d \mu < \infty \} < \infty$. The supremum is equal to $\IP_\Sigma^\ast (\varphi)$ by theorem \ref{nc}.
\ep
\appendix
\section{Pressure at a Point.}
In theorems \ref{ac.3}, \ref{ac.4} and the remark afterwards, we considered the topological pressure on a point $z$. Here, we prove the formulae that we quoted for $P_{ \{z\} } (\varphi )$, $\underline {CP}_{\{z\}} (\varphi)$ and $\overline {CP}_{\{z\}} (\varphi)$. While the proof is an elementary argument direct from from the definitions, it captures a fundamental difference between the different definitions, and it is for this reason that we include it.
\bt \label{app}
Let $X$ be a compact metric space, $f: X \mapsto X$ and $z$ be an arbitrary point. Then
\[
P_{ \{z\} } (\varphi ) = \ul {CP}_{\{z\}} (\varphi) = \liminf_{n \ra \infty} \frac{1}{n} \sum_{i = 0}^{n-1} \varphi (f^i (z)),
\]
\[
\overline {CP}_{\{z\}} (\varphi) = \limsup_{n \ra \infty} \frac{1}{n} \sum_{i = 0}^{n-1} \varphi (f^i (z)).
\]
\et
\begin{remark}
It follows from theorem \ref{app} and the ergodic theorem that for any invariant measure $\mu$, there is a set of full measure so that $\IP_{ \{z\} } (\varphi ) = \ul {CP}_{\{z\}} (\varphi) = \overline {CP}_{\{z\}} (\varphi)$. If $\mu$ is ergodic, this value is $\int \varphi d \mu$. 
\end{remark}
\begin{remark}
If $z$ is a point for which the Birkhoff average of $\varphi$ does not exist, then $\IP_{ \{z\} } (\varphi ) = \ul {CP}_{\{z\}} (\varphi) < \overline {CP}_{\{z\}} (\varphi)$.
\end{remark}
The theorem is a consequence of the lemmas that follow and the relation $\IP_{Z } (\varphi ) \leq \ul {CP}_{Z} (\varphi) \leq \overline {CP}_{Z} (\varphi)$ for any Borel set $Z \subset X$ (formula (11.9) of \cite{Pe}).

\begin{lemma} \label{Ak}
Let $(X,d)$ be a compact metric space, $\varphi: X \mapsto \IR$ a continuous function, and $z \in X$. Then
\[
\IP_{ \{z\} } (\varphi ) \geq \liminf_{n \ra \infty} \frac{1}{n} \sum_{i = 0}^{n-1} \varphi (f^i (z)).
\]
\end{lemma} 
\bp
Let $\mc U_\epsilon$ denote the set of open balls $B(x, \epsilon)$ in $X$. In the definition of $\PZ$, it suffices to consider open covers of the form $\mc U_\epsilon$ and strings of the form $\BU = \{B (x, \epsilon), B(fx, \epsilon), \ldots, B (f^{n-1}x, \epsilon) \}$ (see \cite{Pe}, remark 1 after theorem 11.5). In the notation of the proof of lemma \ref{spec}, $\XU = B_n (x, \epsilon)$. Without loss of generality, it suffices to consider collections $\mc G$ containing only one such string which covers $\{z\}$. We identify such a collection $\mc G$ with a single set $B_n (x, \epsilon)$ which contains $z$. Fix $\epsilon > 0$, $N \in \IN$ and $0 < \delta < \frac{1}{2}$. Choose $\alpha$ satisfying 
\begin{equation} \label{Aj}
\alpha < \liminf_{n \ra \infty} \frac{1}{n} \sum_{i = 0}^{n-1} \varphi (f^i (z)) - \gamma( \epsilon) - \delta,
\end{equation}
where $\gamma( \epsilon) = \sup \{ | \varphi(x) - \varphi (y) | : |x-y| < \epsilon \} $. Assume $N$ was chosen sufficiently large so that for $m \geq N$,
\begin{equation} \label{Ak1}
\frac{1}{n} \sum_{i = 0}^{n-1} \varphi (f^i (z)) \geq \liminf_{n \ra \infty} \frac{1}{n} \sum_{i = 0}^{n-1} \varphi (f^i (z)) - \delta.
\end{equation}
Choose $\mc G = \{ B_m (x, \epsilon) \}$ such that $z \in B_m (x, \epsilon)$, $m \geq N$ and
\[
\left| Q(\{z\},\alpha, \mc U_\epsilon, \mc G, \varphi) - M(\{z\}, \alpha, \mc U_\epsilon, N, \varphi) \right| \leq \delta.
\]
We can prove that $\sum_{k=0}^{m -1} \varphi(f^{k}(z)) - m\gamma( \epsilon) - m\alpha > 0$, which follows from (\ref{Aj}) and (\ref{Ak1}). It follows that
\begin{eqnarray*}
M(\{z\}, \alpha, \mc U_\epsilon, N, \varphi) & \geq &\exp \left \{-\alpha m +\sup_{y\in \mathit{B}_m (x, \epsilon)} \sum_{k=0}^{m -1} \varphi(f^{k}(y))\right \} - \delta \\
& \geq & \exp \left \{-\alpha m + \sum_{k=0}^{m -1} \varphi(f^{k}(z)) - m \gamma( \epsilon) \right \} - \delta \\
& \geq & 1 - \delta \geq \frac{1}{2}.
\end{eqnarray*}
So $m(\{z\}, \alpha, \mc U_\epsilon, \varphi) > 0$ and hence $\IP_Z (\varphi, \mc U_\epsilon) \geq \alpha$. It follows that $\IP_Z (\varphi, \mc U_\epsilon) \geq \liminf_{n \ra \infty} \frac{1}{n} \sum_{i = 0}^{n-1} \varphi (f^i (z)) - \gamma( \epsilon) - \delta$. On taking the limit $\epsilon \ra 0$ and noting that $\delta$ was arbitrary, we obtain the desired result.
\ep
\bl \label{A2}
$\overline {CP}_{\{z\}} (\varphi) =\limsup_{n \ra \infty} \frac{1}{n} \sum_{i = 0}^{n-1} \varphi (f^i (z)).$
\el
\bp
It follows from the definition of $\overline {CP}_{Z} (\varphi)$ that
\[
\ul {CP}_{\{z\}} (\varphi) = \lim_{\epsilon \ra 0} \limsup_{n \ra \infty} \frac{1}{n} \log \left(\inf_{ x : z \in B_n (x, \epsilon) } \exp \left \{\sum_{i = 0}^{n-1} \varphi (f^i (x)) \right\} \right).
\]
For a fixed $\epsilon$ and $B_n (x, \epsilon)$ which contains $z$, 
\[
\sum_{i = 0}^{n-1} \varphi (f^i (x)) \geq \sum_{i = 0}^{n-1} \varphi (f^i (z)) - n \gamma(\epsilon).
\] 
It follows that
\[
\overline {CP}_{\{z\}} (\varphi) \geq \lim_{\epsilon \ra 0} \limsup_{n \ra \infty} \frac{1}{n} \{\sum_{i = 0}^{n-1} \varphi (f^i (z)) - \gamma(\epsilon)\}.
\]
We obtain $\overline {CP}_{\{z\}} (\varphi) \geq \limsup_{n \ra \infty} \frac{1}{n} \sum_{i = 0}^{n-1} \varphi (f^i (z))$. The proof of the reverse inequality is similar.
\ep
\bl
$\ul {CP}_{\{z\}} (\varphi) = \liminf_{n \ra \infty} \frac{1}{n} \sum_{i = 0}^{n-1} \varphi (f^i (z)).$
\el
\bp
It follows from the definition of $\ul {CP}_{Z} (\varphi)$ that
\[
\ul {CP}_{\{z\}} (\varphi) = \lim_{\epsilon \ra 0} \liminf_{n \ra \infty} \frac{1}{n} \log \left(\inf_{ x : z \in B_n (x, \epsilon) } \exp \left \{\sum_{i = 0}^{n-1} \varphi (f^i (x)) \right\} \right).
\]
The rest of the proof proceeds in the same way as that of lemma \ref{A2}.
\ep
\section*{Acknowledgements}
This work was completed at the Mathematics Institute of the University of Warwick and constitutes part of my PhD, which was supported by EPSRC. I would like to thank my supervisors Peter Walters and Mark Pollicott for many useful discussions and reading draft versions of this work, for which I am most grateful. I also wish to thank Yakov Pesin for a helpful discussion in the early stages of this work.

\bibliographystyle{plain}
%\bibliography{pressure_thompson}
\bibliography{master}
\end{document}